\documentclass{article}

\usepackage{amssymb}
\usepackage{amsfonts}
\usepackage{amsmath}
\usepackage{latexsym}

\newtheorem{theorem}{Theorem}[section]

\newtheorem{proposition}[theorem]{Proposition}
\newtheorem{corollary}[theorem]{Corollary}
\newtheorem{lemma}[theorem]{Lemma}
\newtheorem{remark}[theorem]{Remark}
\newtheorem{problem}[theorem]{Problem}
\newtheorem{example}[theorem]{Example}

\newcommand{\FF}{\mathbb F}
\newcommand{\CC}{\mathbb C}
\newcommand{\RR}{\mathbb R}
\newcommand{\NN}{\mathbb N}
\newcommand{\ZZ}{\mathbb Z}

\def\Gl{\mathop{\rm Gl}\nolimits}
\def\deg{\mathop{\rm deg }\nolimits}
\def\rank{\mathop{\rm rank}\nolimits}
\def\diag{\mathop{\rm diag }\nolimits}
\newcommand{\se}{\ensuremath{\stackrel{s.e.}{\sim}}}

\title{Weierstrass structure and eigenvalue placement of regular matrix pencils under low rank perturbations}

\author{Itziar Baraga\~na, 
Departamento de Ciencia de la Computaci\'on e I.A.\\
Facultad de Inform\'atica, Universidad del Pa\'{\i}s Vasco, UPV/EHU
\footnote{itziar.baragana@ehu.eus}\\
\\
Alicia Roca, 
Departamento de  Matem\'atica Aplicada, IMM \\ 
Universitat Polit\`ecnica Val\`encia\\ 46021 Valencia, Spain
\footnote{aroca@mat.upv.es}}
\date{}

\begin{document}
\maketitle

\begin{abstract}
We solve the problem of determining the Weierstrass structure of a regular matrix pencil obtained by a low rank  perturbation of another regular matrix pencil.
We apply the result to find necessary and sufficient conditions for the existence of a low rank perturbation such that the perturbed pencil has prescribed eigenvalues and algebraic multiplicities. The results hold over   fields with sufficient number of elements.
\end{abstract}


AMS:   15A22, 47A55, 15A18

\section{Introduction}
\label{secintroduction}

In the last decades the problem of {\em low rank perturbation} has been 
widely studied by different authors from different points of view. Given a matrix $A$, the problem consists in characterizing the invariants of  $A+P$   with respect  to a given equivalence relation, where  $P$ is a matrix of   bounded rank.
 As  pointed out in \cite{DoSt14}, the problem is equivalent to the {\em rank distance problem}, i.e. the problem of finding  two matrices with prescribed invariants, such that their difference has bounded rank.

Different requirements on the matrices $A, P$ and a matrix $B$ equivalent to $A+P$, on  their domain and on the equivalence relation,   lead to different types of problems.   Additionally, in many contributions in the area the problem is addressed  generically, i.e.  the perturbation $P$ belongs to an open and dense subset of the set of matrices with rank less than or equal to $r$, for a given integer $r$ (see for instance \cite{Batzke14, Batzke15, TeDo07, TeDo16,  TeDoMo08,  MeMeRaRo11, MoDo03, Sa02, Sa04} and the references therein). 
In other cases, the matrix $P$ is an arbitrary perturbation belonging to the whole set of matrices of  rank less than or equal to $r$.  In this paper we follow the second approach. Some  results related to the problem from this point of view can be found  in  \cite{DoSt14, GeTr17, LMPTW18, Sa79, Silva88_2,Silva88_1,Th80, Za91}. We analyze next most of them.

For matrices over a principal ideal domain the problem of characterizing the invariant factors of $A+P$, where $P$ is a matrix with 
$\rank (P)\leq r$, has been solved in \cite{Th80} for $r=1$ and in \cite{Sa79} for the general case (see Lemma \ref{theoSa} below).
For square matrices over a field and  the similarity relation, a  solution is given in \cite{Silva88_1} and \cite{Za91}  (see Proposition \ref{lemmaZa} below). For $r=1$ the problem was already solved in \cite{Th80}.
The case where the perturbation $P$ has fixed rank  ($\rank (P)=r$) has been solved in \cite{Silva88_1}  for square matrices over algebraically closed fields, and for matrices over principal ideal domains in \cite{Silva88_2}. 
More recently, in \cite{DoSt14}, the low rank perturbation problem is solved for pairs of matrices and the feedback equivalence relation.

It is well known that two matrices $A, B\in \FF^{n \times n}$, $\FF$ a field,  are similar if and only if their characteristic matrices $sI_n-A$ and $sI_n-B$ are equivalent as polynomial matrices (see for instance \cite[Ch. 2]{Friedland80}, \cite[Ch. 6]{Ga74}).   Therefore,  the results in \cite{Silva88_1,Th80, Za91} hold for regular pencils having only finite invariant factors,  hence they give a solution to the low rank perturbation problem  for this type of pencils and  for a constant perturbation.
In the same way, given a matrix pair $(A_1,A_2)$ with $A_1\in \FF^{n \times n}$, $A_2 \in \FF^{n \times m}$,  we can associate to $(A_1,A_2)$ the singular matrix pencil $\begin{bmatrix}sI_n-A_1&A_2\end{bmatrix}$. Two matrix pairs are feedback equivalent if and only if the associated pencils are strictly equivalent (see for instance \cite[Ch. IX]{GKvS95}).
Then, the result in \cite{DoSt14} holds for singular matrix pencils of the
form $\begin{bmatrix}sI_n-A_1&A_2\end{bmatrix}$ and the strict equivalence relation of matrix pencils, when the perturbation matrix is constant.

In this paper we  solve the low rank perturbation problem for regular matrix pencils and the strict  equivalence relation, when the perturbation matrix is allowed to  be a pencil ({pencil perturbation problem}). The general pencil perturbation problem for arbitrary singular pencils remains open.

A complete system of invariants for the strict equivalence of regular matrix pencils is formed by the invariant factors (equivalently the finite elementary divisors) and the infinite elementary divisors. We refer to them as the {\em Weierstrass structure} of the pencil.
Given a regular matrix pencil $A(s)$, we obtain necessary conditions for  the Weierstrass structure of the pencil $A(s)+P(s)$, 
 when $P(s)$ is a pencil with $\rank (P(s))\leq r$ and $A(s)+P(s)$ is regular, which hold for arbitrary fields.
Conversely, 
we prove that the  necessary conditions  obtained are sufficient, 
but in this case  we need to impose a condition on the field, we need that the field  has sufficient number of elements (see Remark  \ref{cardF}).

As mentioned above, the solutions to the perturbation problem provided in \cite{Silva88_1,Th80, Za91} are also solutions to our problem in the  case where the pencils do not have  infinite elementary divisors. Notice that  the perturbation matrix P is required there to be constant.
In the more general statement of perturbation of polynomial matrices of \cite{Sa79, Silva88_2, Th80}, the only invariants involved in the equivalence relation are the (finite) invariant factors, and the perturbation matrix can be a polynomial matrix of any degree. 
In this paper, the perturbation problem is solved for arbitrary regular pencils, therefore including the case where they have  infinite elementary divisors, and  the perturbation matrix is allowed to be a polynomial matrix of degree at most one.

Moreover, the result we obtain allows us to generalize  the solution given in \cite{GeTr17}  to the eigenvalue placement problem. Namely, 
 given a regular matrix pencil $A(s)$,
we obtain necessary and sufficient conditions for the existence of a  pencil $P(s)$ with $\rank (P(s))\leq r$ such that $A(s)+P(s)$ is regular with prescribed eigenvalues and algebraic multiplicities. In \cite{GeTr17} the solution to this problem is obtained for $r=1$ and for complex or real pencils.

For complex regular matrix pencils, generic low rank perturbations have been studied for instance in \cite{TeDo16, TeDoMo08}. There, an eigenvalue of  the pencil $A(s)$ is fixed and  the spectral behavior of the eigenvalue after perturbation is generically characterized. In our result  all the possible achievable  partial multiplicities of all of the eigenvalues are characterized, therefore including those of the  generic behavior.
The  same type of problem for  singular matrix pencils  was studied in \cite{TeDo07} and   for  structured regular pencils, for example,  in \cite{Batzke14, Batzke15}. In both cases, the problem is addressed generically.

The paper is organized as follows. In Section \ref{secpreliminaries} we introduce the notation, basic definitions and preliminary results. In Section \ref{secproblem} the problems  to be studied in the paper are established. In Section \ref{secmain} we provide in  Theorem \ref{mainth} a solution to the stated low rank perturbation problem. Section \ref{secplacement} is devoted to solve the eigenvalue placement problem, and the solution is given in Theorem \ref{mainthdet}. Finally, in Section \ref{secconclusions} we summarize the main contributions of the paper.

\section{Preliminaries} 
\label{secpreliminaries}
Let $\FF$ be a field. $\FF[s]$ denotes the ring of polynomials in the indeterminate $s$ with coefficients in $\FF$ and  $\FF[s, t]$  the ring of polynomials in two
variables $s, t$ with coefficients in $\FF$.
We denote by $\FF^{m\times n}$,  $\FF[s]^{m\times n}$ and $\FF[s, t]^{m\times n}$ the vector spaces  of $m\times n$ matrices with elements in $\FF$, $\FF[s]$ and $\FF[s, t]$, respectively.
$\Gl_n(\FF)$ will be the general linear group of invertible matrices
in $\FF^{n\times n}$. 
A matrix    $U(s)\in \FF[s]^{n\times n}$ is  {\em unimodular} if $0\neq \det (U(s))\in \FF$, i.e. it is a unit in the ring  $\FF[s]^{n\times n}$.

\medskip
Given a polynomial matrix $G(s)\in \FF[s]^{m\times n}$, the {\em degree} of $G(s)$, denoted by $\deg(G(s))$, is the  maximum  of  the degrees of its entries. The {\em normal rank} of $G(s)$, denoted by $\rank (G(s))$,  is the order of the largest non identically zero minor of $G(s)$, i.e. it is the rank of $G(s)$ considered as a matrix on the field of fractions of $\FF[s]$. The {\em determinantal divisor of order $k$} of  $G(s)$, denoted by $D_k(s)$, is the monic greatest common divisor  of the minors of order  $k$ of $G(s)$,  $1\leq k \leq \rank (G(s))$.

\medskip
Two polynomial matrices $G(s), H(s)\in \FF[s]^{m\times n}$ are {\em equivalent} ($G(s)\sim H(s)$) if there exist unimodular matrices $U(s)\in \FF[s]^{m\times m}$, $ V(s) \in \FF[s]^{n\times n}$ such that 
$G(s)=U(s)H(s)V(s)$. 
If $G(s)\in \FF[s]^{m\times n}$ and
$\rank (G(s))=\rho$,
it is well known (see for example \cite[Ch. 6]{Ga74}) that $G(s)$ is equivalent to a unique matrix of the form
$$S(s)=\begin{bmatrix}\diag(\gamma_1(s), \dots, \gamma_\rho(s))&0\\0&0\end{bmatrix},$$ 
where $\gamma_1(s), \dots,\gamma_\rho(s)$ are monic polynomials and
$\gamma_1(s)\mid \dots\mid \gamma_\rho(s)$.
Moreover, 
$$
D_k(s)= \gamma_1(s) \dots\gamma_k(s), \quad 1\leq k \leq \rho,
$$
which means that
$$
\gamma_k(s)=\frac{D_k(s)}{D_{k-1}(s)}, \quad 1\leq k \leq \rho \quad (D_0(s)=1).
$$
The matrix $S(s)$ is the  {\em Smith form} of $G(s)$ and the polynomials $\gamma_1(s), \dots,  \gamma_\rho(s)$ are the  {\em invariant factors} of $G(s)$.
We will take $\gamma_i(s):=1$ for $i<1$ and  $\gamma_i(s):=0$ for $i>\rho$. 

The invariant factors form a complete system of invariants for the equivalence of polynomial matrices, i.e. two polynomial matrices $G(s), H(s)\in \FF[s]^{m\times n}$ are equivalent if and only if they have the same invariant factors.

\medskip
A {\em matrix pencil} is a  polynomial matrix $G(s)\in \FF[s]^{m\times n}$
such that $\deg(G(s))\leq1$. The pencil is  {\em regular} if $m=n$ and $\det(G(s))$ is not the zero polynomial. Otherwise it is  {\em singular}.

Two matrix pencils 
$G(s)=G_0+sG_1, H(s)=H_0+sH_1\in \FF[s]^{m\times n}$ are {\em strictly equivalent} ($G(s)\se H(s)$) 
if there exist invertible matrices $Q\in \Gl_m(\FF)$,   $R\in \Gl_n(\FF)$ such that 
$G(s)=QH(s)R$.
Equivalently, 
$
G_0=QH_0R,  
$
$
G_1=QH_1R
$.

It is immediate that if $G(s)\se H(s)$ then $G(s)\sim H(s)$.
Moreover, if $n=m$, $\det(G_1)\neq 0$ and $\det(H_1)\neq 0$, then $G(s)\se H(s)$ if and only if
$G(s)\sim H(s)$
(see for instance \cite[Ch. 12, Theorem 1]{Ga74}). 

\medskip
If $\gamma_1(s)\mid \dots \mid \gamma_\rho(s)$ are the invariant factors of $G(s)=G_0+sG_1\in \FF[s]^{m\times n}$,
then the invariant  factors of  the
matrix pencil $\bar G(t)=tG_0+G_1\in \FF[t]^{m\times n}$ can be expressed as
$$
\bar \gamma_i(t)=k_it^{q_i}t^{\deg(\gamma_i)}\gamma_i(\frac{1}{t}),\quad 1\leq i \leq \rho,
$$
for some integers $q_i\geq 0$, where $0\neq k_i\in \FF$ are scalars such that $\bar \gamma_i(t)$ are monic. If  $q_i>0$, then  $t^{q_i}$ is an {\em infinite elementary divisor} of $G(s)$.
The infinite elementary divisors of $G(s)$  exist if and only if $\rank(G_1)< \rank (G(s))$.

The exponents $q_i$ are the {\em partial multiplicities at infinity of $G(s)$}
and we  will denote them by $m_i(\infty, G(s)):=q_i$.
Then, $m_1(\infty, G(s))\leq \dots\leq  m_\rho(\infty, G(s))$.

\medskip
Given  $G(s)=G_0+sG_1\in \FF[s]^{m\times n}$,  with $\rank( G(s))=\rho$, the {\em homogeneous pencil associated to $G(s)$} is
$$G(s, t)=tG_0+sG_1\in \FF[s, t]^{m\times n},$$
and the {\em homogeneous determinantal divisor of order $k$} of  $G(s)$, denoted by $\Delta_k(s, t)$, is the  greatest common divisor  of the minors of order  $k$ of $G(s, t)$,
$1\leq k \leq \rho$. We will assume that $\Delta_k(s, t)$ is monic with respect to $s$.
The homogeneous determinantal divisors of $G(s)$ are homogeneous polynomials and $\Delta_{k-1}(s, t)\mid \Delta_k(s,t)$, $1\leq k \leq \rho$.
Defining
$$
\Gamma_k(s, t)=\frac{\Delta_k(s,t)}{\Delta_{k-1}(s, t)}, \quad 1\leq k \leq \rho \quad (\Delta_0(s,t)=1),
$$
if $\gamma_1(s)\mid \dots \mid \gamma_\rho(s)$ and $\bar \gamma_1(t)\mid \dots \mid \bar \gamma_\rho(t)$ are the invariant factors of $G(s)$ and  
$\bar G(t)=tG_0+G_1$, respectively, then
$$
\gamma_i(s)=\Gamma_i(s,1), \quad \bar \gamma_i(t)=k_i\Gamma_i(1,t), \quad  1\leq i \leq \rho \quad (0\neq k_i\in \FF),
$$
and
$$
\Gamma_i(s,t)=t^{m_i(\infty, G(s))}t^{\deg(\gamma_i)}\gamma_i(\frac{s}{t}),\quad 1\leq i \leq \rho.
$$ 
As a consequence, $\Gamma_1(s,t)\mid \dots \mid\Gamma_{\rho}(s, t)$.
The polynomials $\Gamma_1(s,t), \dots, \Gamma_{\rho}(s, t)$ are called
the {\em homogeneous invariant factors of $G(s)$}. For details see  \cite[Ch. 2]{Friedland80}, \cite[Ch. 12]{Ga74}.
We will take $\Gamma_i(s,t):=1$ for $i<1$ and 
$\Gamma_i(s,t):=0$ for $i>\rho$. 
Observe that the homogeneous invariant factors of  $\bar G(t)=tG_0+G_1$ are
$\bar \Gamma_i(s,t)=\Gamma_i(s,t)$.

The homogeneous invariant  factors  form a complete system of invariants for the strict equivalence of regular pencils. A proof of the following theorem can be found,  for instance,  in
\cite[Ch. 12]{Ga74} for infinite fields and in \cite[Ch. 2]{Ro03} for arbitrary fields.
\begin{theorem}[\mbox{Weierstrass}]
  Two regular matrix pencils are strictly equivalent if and only if they have the same homogeneous invariant factors.
\end{theorem}
Notice that the invariant factors and the infinite elementary divisors determine the homogeneous invariant factors. As a consequence, two regular matrix pencils are strictly equivalent if and only if they have the same  invariant factors and the same infinite elementary divisors.

\medskip

We denote by $\overline{ \FF}$ the algebraic clousure of $\FF$. 
The {\em finite spectrum} of a regular pencil
$G(s)=G_0+sG_1\in \FF[s]^{n\times n}$ is defined as
$$
\Lambda_f(G(s))=\{\lambda\in \overline{\FF}\; : \; \det(G(\lambda))=0\}.
$$
If $\gamma_1(s)\mid \dots\mid \gamma_n(s)$ are the invariant factors of  $G(s)$,
then we can write
$$
\gamma_i(s)=\prod_{\lambda\in \Lambda_f(G(s))}(s-\lambda)^{m_i(\lambda, G(s))}, 
$$
where 
$0\leq m_1(\lambda, G(s))\leq \dots\leq   m_n(\lambda, G(s))$ are called the 
{\em partial multiplicities at $\lambda$ of $G(s)$}.
The {\em spectrum} of 
$G(s)$ is $\Lambda(G(s))=\Lambda_f(G(s))$ if $\det(G_1)\neq 0$ and
$\Lambda(G(s))=\Lambda_f(G(s))\cup\{\infty\}$ if $\det(G_1)= 0$.
The elements $\lambda\in \Lambda(G(s))$ are the {\em eigenvalues} of $G(s)$.

If  $\lambda \in \overline{\FF}\setminus \Lambda(G(s))$, 
we put $m_1(\lambda, G(s))=\dots=m_n(\lambda,  G(s))=0$.
For $\lambda \in \overline{\FF}\cup \{\infty\}$, we  will agree that $m_i(\lambda, G(s))=0$ for $i<1$ and  $m_i(\lambda, G(s))=\infty$ for $i>n$. 
We denote by
$(w_1(\lambda, G(s)), \dots, w_n(\lambda, G(s))$ the conjugate partition of
the partition
$(m_n(\lambda, G(s)), \dots, m_1(\lambda, G(s))$, i.e.
$$
w_i(\lambda, G(s))=\#\{j\in\{1, \dots, n\}\; : \; m_j(\lambda, G(s))\geq i\}, \quad
1\leq i \leq n.
$$
The {\em algebraic multiplicity} of $\lambda$ in  $G(s)$ is  $\mu_a(\lambda, G(s)):= \sum _{i=1}^n m_i(\lambda, G(s))$, and the {\em geometric multiplicity} is $\mu_g(\lambda, G(s)):=\#\{i\in\{1, \ldots, n\}: m_i(\lambda, G(s))>0\}$, i.e. 
$\mu_g(\lambda, G(s))=w_1(\lambda, G(s))$.

The following technical result characterizes some inequality relations between  the elements of two partitions of nonnegative integers and those of its conjugate partitions. 
\begin{lemma}\cite[Lemma 3.2]{LiSt09}\label{lemmaconj}
Let 
$p_1\geq p_2\geq \dots \geq 0$
and
$p'_1\geq p'_2\geq \dots \geq 0$
be partitions of $n$ and $n'$
with conjugate partitions
$q_1\geq q_2\geq \dots\geq 0$
and
$q'_1\geq q'_2\geq \dots \geq 0$. Let $r\in \NN$. Then 
$q'_i\geq q_{i+r}$ and $q_i\geq q'_{i+r}$ for all $i>0$ if and only if $\mid p_i-p'_i\mid \leq r$ for all $i>0$.
\end{lemma}

\medskip
In the next lemma, whose proof is straightforward, we show that 
 conditions of divisibility between homogeneous invariant factors can be expressed in terms of 
 divisibility between invariant factors and infinite elementary divisors.
\begin{lemma}\label{lemmadiv}
Let $\gamma(s),\omega(s)\in \FF[s]$ be monic polynomials, and let $m, m'$ be nonnegative integers.
If
$\Gamma(s,t)=t^{m}t^{\deg(\gamma)}\gamma(\frac{s}{t})$ and 
$\Omega(s,t)=t^{m'}t^{\deg(\omega)}\omega(\frac{s}{t})$, then
$$\Gamma(s,t)\mid \Omega(s,t) \, \mbox{ 
if and only if }\, 
\left\{\begin{array}{l}
\gamma(s)\mid \omega(s), \\
m\leq m'.
\end{array}\right.
$$
\end{lemma}

\section{Statement of the problems}
\label{secproblem}

The first problem we deal with in this paper is the following one.

\begin{problem}[Low rank perturbation for regular matrix pencils]\label{problem}
Given two regular matrix pencils  $A(s), B(s)\in \FF[s]^{n \times n}$
and a nonnegative integer $r$,
 find necessary and sufficient conditions for the existence of a 
matrix pencil
 $P(s)\in \FF[s]^{n \times n}$ such that $\rank( P(s))\leq r$ and $A(s)+P(s)\se B(s)$.  
\end{problem}

Next lemma shows that the pencil  $A(s)$ can  be substituted by any other 
pencil
strictly equivalent to $A(s)$.

\begin{lemma}\label{lemmasust}
 Let  $A(s), B(s), P(s)\in \FF[s]^{n \times n}$ be matrix pencils.
Let $Q, R\in \Gl_n(\FF)$ and $ A'(s)=Q A(s)R$. If   $A(s)+P(s)\se B(s)$ then $A'(s)+QP(s)R\se B(s)$.
\end{lemma}

    {\bf Proof.}
If $B(s)\se A(s)+P(s)$, then $B(s)\se Q(A(s)+P(s))R=A'(s)+QP(s)R$. 

\hfill $\Box$
  
In the next lemma we see that the roles of the pencils $A(s)$ and $B(s)$ can be interchanged.
\begin{lemma}\label{lemmasimetry}
Let $A(s), B(s), P(s)\in \FF[s]^{n \times n}$ be matrix pencils such that  $A(s)+P(s)\se B(s)$. Then there exists a matrix pencil $P'(s)\in \FF[s]^{n \times n}$ such that  $B(s)+P'(s)\se A(s)$ and $\rank( P'(s))=\rank (P(s))$.
\end{lemma}

{\bf Proof.}
There exist $Q, R\in \Gl_n(\FF)$ such that 
$Q(A(s)+P(s))R=B(s)$. Hence,
$B(s)-QP(s)R=QA(s)R\se A(s)$ and the lemma follows with $P'(s)=-QP(s)R$. 

\hfill $\Box$

The solution to Problem \ref{problem} will allow us to also solve the following one. 

\begin{problem}[Eigenvalue placement for regular matrix pencils under low rank perturbations]\label{prdet}
  Given a regular matrix pencil  $A(s)\in \FF[s]^{n \times n}$, a nonnegative integer $r$ and a monic polynomial  $0\neq p(s)\in \FF[s]$ with $\deg(p(s))\leq n$,  
 find necessary and sufficient conditions for the existence of a 
 matrix pencil
 $P(s)\in \FF[s]^{n \times n}$ such that $\rank (P(s))\leq r$ and   
$\det(A(s)+P(s))=kp(s)$, with $k\in \FF$.
\end{problem}

Notice that in this problem  we prescribe the spectrum $\Lambda(A(s)+P(s))$
and the algebraic multiplicities of every $\lambda \in \Lambda(A(s)+P(s))$.
For $r=1$ and $\FF=\CC$ or $\FF=\RR$, a solution to Problem \ref{prdet} is given in  \cite{GeTr17}.

\section{Low rank perturbation for regular matrix pencils}\label{secmain}

In this section we give a solution to Problem \ref{problem}.

\subsection{Necessary conditions}\label{subsecnec}

The next lemma was obtained in \cite{Sa79}. The result can also be found in \cite[Theorem 1]{Silva88_2}.

\begin{lemma} 
[\mbox{\cite[Theorem 6.1]{Sa79}}]\label{theoSa}
Let $G(s), H(s)\in \FF[s]^{m\times n}$, $n\leq m$,  be polynomial matrices with
invariant factors 
$\gamma_1(s)\mid \dots \mid \gamma_n(s)$ 
and $\delta_1(s)\mid \dots \mid \delta_n(s)$, respectively.
Let $r$ be a nonnegative integer.
There exist matrices $\hat{G}(s), \hat{H}(s)\in \FF[s]^{m\times n}$ equivalent to $G(s)$ and $H(s)$, respectively, such that 
$$
\rank(\hat{H}(s)-\hat{G}(s))\leq r
$$
 if and only if
\begin{equation}\label{eqSa}
\delta_i(s)\mid \gamma_{i+r}(s), \quad \gamma_i(s)\mid \delta_{i+r}(s), \quad 1\leq i \leq n-r.
\end{equation}
  \end{lemma}

\begin{remark}\label{remSa}
Condition (\ref{eqSa})
is equivalent to
$$
 \gamma_{i-r}(s)\mid \delta_i(s)\mid\gamma_{i+r}(s), \quad 1\leq i \leq n,
$$
and to
$$
\delta_{i-r}(s)\mid \gamma_i(s)\mid\delta_{i+r}(s), \quad 1\leq i \leq n.
$$
\end{remark}

Taking advantage of the above result,   we obtain  desired necessary  conditions for solving  Problem \ref{problem} in the next  proposition. 

\begin{proposition}\label{propnecesr}
Let $A(s)=A_0+sA_1, B(s)=B_0+sB_1\in \FF[s]^{n \times n}$ be regular matrix pencils. Let $P(s)=P_0+sP_1\in \FF[s]^{n \times n}$ be a matrix pencil with  $ \rank (P(s))= r$ such that  $A(s)+P(s)\se B(s)$.
Let $\phi_1(s,t)\mid \dots\mid \phi_n(s, t)$ and $\psi_1(s,t)\mid \dots\mid \psi_n(s, t)$ be the homogeneous invariant factors of 
$A(s)$ and $B(s)$, respectively. Then
\begin{equation}\label{interlacinghomogr1}
\phi_{i-r}(s, t)\mid \psi_i(s, t)\mid\phi_{i+r}(s, t), \quad 1\leq i \leq n.
\end{equation}

\end{proposition}
{\bf Proof.}
Let $\bar A(t)=tA_0+A_1$,  $\bar B(t)=tB_0+B_1$, $\bar P(t)=tP_0+P_1$.

Since $A(s)+P(s)\se B(s)$, we have that  $A(s)+P(s)\sim B(s)$ and $\bar A(t)+\bar P(t)\sim \bar B(t)$.

Let $\alpha_1(s)\mid \dots\mid \alpha_n(s)$ and $\beta_1(s) \mid \dots\mid \beta_n(s)$ be 
the invariant factors of $A(s)$ and $B(s)$, respectively.
The invariant factors of $\bar A(t)=tA_0+A_1$ and $\bar B(t)=tB_0+B_1$ are, respectively,
$$
\bar \alpha_i(t)=k_it^{m_i(\infty, A(s))}t^{\deg(\alpha_i)}\alpha_i(\frac{1}{t}),\quad 1\leq i \leq n, \quad (0\neq k_i\in \FF),
$$
$$
\bar \beta_i(t)=k'_it^{m_i(\infty, B(s))}t^{\deg(\beta_i)}\beta_i(\frac{1}{t}),\quad 1\leq i \leq n, \quad (0\neq k'_i\in \FF).
$$
By Lemma \ref{theoSa}, 
\begin{equation}\label{eqinalpha}
\alpha_{i-r}(s)\mid \beta_i(s)\mid\alpha_{i+r}(s), \quad 1\leq i \leq n,
\end{equation}
and
\begin{equation}\label{eqinbeta}
\bar \alpha_{i-r}(t)\mid \bar \beta_i(t)\mid \bar \alpha_{i+r}(t),\quad 1\leq i \leq n.
\end{equation}
It follows from (\ref{eqinbeta}) that
\begin{equation}\label{eqinm}
 m_{i-r}(\infty, A(s))\leq m_i(\infty, B(s))\leq m_{i+r}(\infty, A(s)),\quad 1\leq i \leq n.
\end{equation}
By Lemma \ref{lemmadiv}, conditions (\ref{eqinalpha}) and (\ref{eqinm}) are equivalent to 
(\ref{interlacinghomogr1}).

\hfill $\Box$

\begin{remark}\label{r1tor}
If $ \rank (P(s))= r_1\leq r$, this proposition tells us that 
$$\phi_{i-r_1}(s, t)\mid \psi_i(s, t)\mid\phi_{i+r_1}(s, t), \quad 1\leq i \leq n,$$
therefore
$$\phi_{i-r}(s, t)\mid\phi_{i-r_1}(s, t)\mid \psi_i(s, t)\mid\phi_{i+r_1}(s, t)\mid \phi_{i+r}(s, t), \quad 1\leq i \leq n,$$
hence, condition (\ref{interlacinghomogr1}) is necessary for Problem \ref{problem}.
\end{remark}

\begin{remark}\label{remsimhom}
Condition (\ref{interlacinghomogr1})
is also equivalent to
$$
 \psi_{i-r}(s,t)\mid \phi_i(s,t)\mid\psi_{i+r}(s,t), \quad 1\leq i \leq n.
$$

\end{remark}

Condition (\ref{interlacinghomogr1}) can be stated in terms of the partial multiplicities of the elements of $ \Lambda(A(s))\cup \Lambda(B(s))$.

\begin{corollary}\label{cornecmult}
Let $A(s), B(s)\in \FF[s]^{n \times n}$ be regular matrix pencils. Let $P(s)\in \FF[s]^{n \times n}$ be a matrix pencil with  $ \rank (P(s))= r$ such that  $A(s)+P(s)\se B(s)$.
Then
\begin{equation}\label{eqnecmult}
m_{i-r}(\lambda, A(s))\leq m_i(\lambda, B(s))\leq m_{i+r}(\lambda, A(s)), 
\, 1\leq i \leq n,  \quad \lambda \in \overline{\FF}\cup\{\infty\}.
\end{equation}
Equivalently, 
$$m_{i-r}(\lambda, B(s))\leq m_i(\lambda, A(s))\leq m_{i+r}(\lambda, B(s)), 
\, 1\leq i \leq n, \quad  \lambda \in \overline{\FF}\cup\{\infty\}.$$
\end{corollary}

From Lemma \ref{lemmaconj}  we conclude  the following result
(Corollary \ref{coreqinw} for $r=1$ is proved in Proposition 4.2 of \cite{GeTr17}).
\begin{corollary} \label{coreqinw}
Let $A(s), B(s)\in \FF[s]^{n \times n}$ be regular matrix pencils. Let $P(s)\in \FF[s]^{n \times n}$ be a matrix pencil with  $ \rank( P(s))= r$ such that  $A(s)+P(s)\se B(s)$.
Then
\begin{equation}\label{eqinw}
w_{i}(\lambda, A(s))-r\leq w_i(\lambda, B(s))\leq w_i(\lambda, A(s))+r,
\quad 1\leq i \leq n, \quad \lambda \in \overline{\FF}\cup\{\infty\}.
\end{equation}
Equivalently, 
$$
w_{i}(\lambda, B(s))-r\leq w_i(\lambda, A(s))\leq w_i(\lambda, B(s))+r,
\quad 1\leq i \leq n, \quad \lambda \in \overline{\FF}\cup\{\infty\}.
$$
\end{corollary}

\begin{corollary}
Let $A(s), B(s)\in \FF[s]^{n \times n}$ be regular matrix pencils. Let $P(s)\in \FF[s]^{n \times n}$ be a matrix pencil with  $ \rank ( P(s))= r$ such that  $A(s)+P(s)\se B(s)$.
Then
$$
\mu_g(\lambda, A(s))-r\leq \mu_g(\lambda, B(s))\leq \mu_g(\lambda, A(s))+r, \quad \lambda \in \overline{\FF}\cup\{\infty\}.
$$

\end{corollary}

\subsection{Sufficiency of the conditions}\label{subsecsuf}

For the case when $A(s)=sI_n-A$, $B(s)=sI_n-B$, with $A, B\in \FF^{n \times n}$, a solution to Problem \ref{problem} is given in \cite{Th80} for $r=1$ and in
\cite{Silva88_1} and \cite{Za91} for the general case. Taking advantage of this result (see Proposition \ref{lemmaZa}) we  prove in Corollary \ref{case1} the sufficiency of  conditions (\ref{interlacinghomogr1}) for pencils without infinite elementary divisors. 

\begin{proposition} 
[\mbox{\cite[Theorem 1]{Silva88_1}, \cite[Theorem 3]{Za91}}]\label{lemmaZa}
Let $G\in \FF^{n\times n}$ and $\gamma_1(s)\mid \dots \mid \gamma_n(s)$ be its 
invariant factors. Let $\delta_1(s)\mid \dots \mid \delta_n(s)$ be monic polynomials such that $\sum_{i=1}^n\deg(\delta_i(s))=n$. Let $r$ be a nonnegative integer.
Then there exists a matrix
$P\in \FF^{n\times n}$ such that  $\rank( P)\leq r$  and $G+P$ has $\delta_1(s)\mid \dots \mid \delta_n(s)$ as invariant factors if and only if
$$
\delta_{i-r}(s)\mid \gamma_i(s)\mid\delta_{i+r}(s), \quad 1\leq i \leq n.
$$
\end{proposition}
\begin{corollary}\label{case1}
Let $A(s)=A_0+sA_1, B(s)=B_0+sB_1\in \FF[s]^{n \times n}$ be such that 
  $\det(A_1)\neq 0$ and  $\det(B_1)\neq 0$.
Let $\phi_1(s,t)\mid \dots\mid \phi_n(s, t)$ and $\psi_1(s,t)\mid \dots\mid \psi_n(s, t)$ be the homogeneous invariant factors of $A(s)$ and $B(s)$, respectively. Let $r$ be a nonnegative integer. If (\ref{interlacinghomogr1}) is satisfied, then there exists $P_0\in \FF^{n \times n}$ such that $\rank (P_0)\leq r$ and 
$A(s)+P_0\se B(s)$.
\end{corollary}

{\bf Proof.}
Recall that two square matrices $G, H\in \FF^{n \times n}$ are similar if and only if $(sI_n-G)\sim (sI_n-H)$ and that the invariant factors of $G\in \FF^{n \times n}$ are those of $sI_n-G$.

Let $\alpha_1(s)\mid \dots\mid \alpha_n(s)$ and $\beta_1(s), \mid \dots\mid \beta_n(s)$ be 
the invariant factors of $A(s)$  and $B(s)$, respectively.
Since $\det(A_1)\neq 0$ and  $\det(B_1)\neq 0$, these pencils do not have infinite elementary divisors, therefore
$$
m_1(\infty, A(s))=\dots =m_n(\infty, A(s))=0, \quad m_1(\infty, B(s))=\dots =m_n(\infty, B(s))=0.
$$
As a consequence, (\ref{interlacinghomogr1}) is equivalent to 
\begin{equation}\label{eqinalphar}
\alpha_{i-r}(s)\mid \beta_i(s)\mid\alpha_{i+r}(s), \quad 1\leq i \leq n.
\end{equation}
We have that $A(s)\se A_1^{-1}(A_0+A_1s)= A_1^{-1}A_0+sI_n$ and $B(s)\se B_1^{-1}(B_0+B_1s)= B_1^{-1}B_0+sI_n$, hence the invariant factors of $-A_1^{-1}A_0$ and $-B_1^{-1}B_0$ are $\alpha_1(s)\mid \dots\mid \alpha_n(s)$ and $\beta_1(s), \mid \dots\mid \beta_n(s)$, respectively.
 By Proposition \ref{lemmaZa}, there exists $P\in \FF^{n \times n}$ such that $\rank (P)\leq r$ and 
$P+A_1^{-1}A_0+sI_n\sim  B_1^{-1}B_0+sI_n$, hence $P+A_1^{-1}A_0+sI_n\se B(s)$.
Setting $P_0=A_1 P$, we have that $\rank (P_0)\leq r$ and 
by Lemma \ref{lemmasust}, $A(s)+P_0\se B(s)$.

\hfill $\Box$

As an immediate consequence of Corollary \ref{case1}, we obtain in Corollary \ref{case2}
 the sufficiency of conditions (\ref{interlacinghomogr1}) for pencils not having the eigenvalue zero.

\begin{corollary}\label{case2}
Let $A(s)=A_0+sA_1, B(s)=B_0+sB_1\in \FF[s]^{n \times n}$ be  such that 
  $\det(A_0)\neq 0$ and  $\det(B_0)\neq 0$.
Let $\phi_1(s,t)\mid \dots\mid \phi_n(s, t)$ and $\psi_1(s,t)\mid \dots\mid \psi_n(s, t)$ be the homogeneous invariant factors of  $A(s)$ and $B(s)$, respectively. Let $r$ be a nonnegative integer. If (\ref{interlacinghomogr1}) is satisfied, then there exists $P_1\in \FF^{n \times n}$ such that $\rank (P_1)\leq r$ and 
$A(s)+sP_1\se B(s)$.
\end{corollary}

{\bf Proof.}
The homogeneous invariant factors of the pencils $\bar A(t)=tA_0+A_1$
and $\bar B(t)=tB_0+B_1$ are
$\phi_1(s, t), \dots, \phi_n(s, t)$ and
$\psi_1(s, t), \dots, \psi_n(s, t)$, respectively.

By Corollary \ref{case1}, there exists $P_1\in \FF^{n \times n}$ such that $\rank (P_1)\leq r$ and $\bar A(t)+P_1=A_1+P_1+tA_0\se \bar B(t)=B_1+tB_0$.
Therefore $ A(s)+sP_1=s(A_1+P_1)+A_0\se sB_1+B_0=B(s)$.

\hfill $\Box$

The restriction on the field $\FF$ introduced in the next corollary allows us to  perform a change of variable over an arbitrary regular pencil, which turns it into a new pencil without the  zero eigenvalue.

\begin{corollary}\label{case3}
Let  $A(s)=A_0+sA_1, B(s)=B_0+sB_1\in \FF[s]^{n \times n}$ be regular pencils, and assume that  there exists $c\in \FF$ such that $c \not \in \Lambda_f(A(s))\cup \Lambda_f(B(s))$.
Let $\phi_1(s,t)\mid \dots\mid \phi_n(s, t)$ and $\psi_1(s,t)\mid \dots\mid \psi_n(s, t)$ be the homogeneous invariant factors of  $A(s)$ and $B(s)$, respectively. Let $r$ be a nonnegative integer. If  (\ref{interlacinghomogr1}) is satisfied, then there exists $P\in \FF^{n \times n}$ such that $\rank (P)\leq r$ and  $A(s)+(s-c)P\se B(s)$.
\end{corollary}

{\bf Proof.}
 As $c\not \in \Lambda_f(A(s))$ and $c\not \in \Lambda_f(B(s))$, we have that $\det(A_0+cA_1)\neq 0$
and $\det(B_0+cB_1)\neq 0$.

Let us  consider the regular matrix pencils
$$
A'(y):=(A_0+cA_1)+yA_1=A_0+(y+c)A_1=A(y+c)\in \FF[y]^{n \times n},
$$
$$
B'(y):=(B_0+cB_1)+yB_1=B_0+(y+c)B_1=B(y+c)\in \FF[y]^{n \times n}.
$$
The  associated homogeneous pencils are
$$A'(y, z)=z(A_0+cA_1)+yA_1=zA_0+(cz+y)A_1
\in \FF[z,y]^{n \times n},$$
$$B'(y, z)=
z(B_0+cB_1)+yB_1=zB_0+(cz+y)B_1\in \FF[z,y]^{n \times n}.$$
Therefore, the homogeneous invariant factors of $A'(y)$ and $B'(y)$ are
$$
\phi'_i(y, z)=\phi_i(cz+y, z),
\quad 1\leq i \leq n,
$$
and
$$
\psi'_i(y, z)=\psi_i(cz+y, z), \quad 1\leq i \leq n,
$$
respectively. Then, condition (\ref{interlacinghomogr1}) implies
$$\phi'_{i-r}(y, z)\mid \psi'_i(y, z)\mid\phi'_{i+r}(y, z), \quad 1\leq i \leq n.$$
By Corollary \ref{case2}, there exists $P\in \FF^{n \times n}$ such that $\rank (P)\leq r$ and $A'(y)+yP\se B'(y)$.
Then $A(s)+(s-c)P=A'(s-c)+(s-c)P\se B'(s-c)=B(s)$.

\hfill $\Box$

Next theorem gives a  complete solution to Problem \ref{problem} under a restriction on the field $\FF$.
\begin{theorem}\label{mainth}
 Let $A(s), B(s)\in \FF[s]^{n \times n}$ be regular matrix pencils.
Let $\phi_1(s,t)\mid \dots\mid \phi_n(s, t)$ and $\psi_1(s,t)\mid \dots\mid \psi_n(s, t)$ be the homogeneous invariant factors of 
$A(s)$ and $B(s)$, respectively, and assume that  $ \FF\cup \{\infty\}\not \subseteq \Lambda(A(s))\cup \Lambda(B(s))  $. Let $r$ be a nonnegative integer. 
There exists a matrix pencil  $P(s)\in \FF[s]^{n \times n}$ such that $\rank (P(s))\leq r$ and $A(s)+P(s)\se B(s)$ if and only if
(\ref{interlacinghomogr1}) holds.

\end{theorem}
{\bf Proof.}
The necessity follows from Proposition \ref{propnecesr} and Remark \ref{r1tor}.
Assume that (\ref{interlacinghomogr1}) holds.
As $ \FF\cup \{\infty\}\not \subseteq \Lambda(A(s))\cup \Lambda(B(s)) $, there exists $c\in  \FF\cup \{\infty\}$ such that
$c\not \in \Lambda(A(s))\cup \Lambda(B(s))$.
If $c=\infty$, then we apply Corollary \ref{case1} and, if $c\neq \infty$, then we apply Corollary \ref{case3}.

\hfill $\Box$

In Corollary \ref{cormainmult}  we restate Theorem \ref{mainth} in terms of the partial multiplicities of the elements of $ \Lambda(A(s))\cup \Lambda(B(s))$.

\begin{corollary}\label{cormainmult}
Let $A(s), B(s)\in \FF[s]^{n \times n}$ be regular matrix pencils.
Assume that  $ \FF\cup \{\infty\}\not \subseteq \Lambda(A(s))\cup \Lambda(B(s))  $.
Let $r$ be a nonnegative integer.
There exists a matrix pencil  $P(s)\in \FF[s]^{n \times n}$ such that $\rank (P(s))\leq r$ and $A(s)+P(s)\se B(s)$ if and only if
(\ref{eqnecmult}) (equivalently (\ref{eqinw}))  holds.

\end{corollary}

\begin{remark}\label{cardF}

In Theorem \ref{mainth}, the condition on the field $\FF$ means that
there exists an element $c\in \FF\cup \{\infty\}$ which is neither an eigenvalue of  $A(s)$ nor of $B(s)$. Observe that if $\#\FF> 2n $, this condition is automatically satisfied. In  the case that    $\#\FF\leq 2n $, Theorem \ref{mainth} can still be applied, as we show   in Example \ref{firstexample}.
Moreover,  the condition $ \FF\cup \{\infty\} \not \subseteq \Lambda(A(s))\cup \Lambda(B(s))$ is not always necessary. 
In Theorem \ref{teoeqmult} below  we prove that even in the case that  $\FF\cup \{\infty\} \subseteq \Lambda(A(s))\cup \Lambda(B(s)) $, 
for certain special structures of $A(s)$ and $B(s)$, condition (\ref{interlacinghomogr1}) is sufficient. We illustrate the case in Example \ref{secondexample}.
\end{remark}

\begin{example}\label{firstexample}
This example is valid for any field $\FF$, including $\FF=\ZZ_2$.

Let
$$
A(s)=\begin{bmatrix}s&0&0\\0&s&1\\0&0&s\end{bmatrix},\quad
B(s)=\begin{bmatrix}1&s&0\\0&1&0\\0&0&s\end{bmatrix}.
$$
The homogeneous invariant factors of $A(s)$ and $B(s)$ are
$
\phi_1(s,t)=1,  \phi_2(s,t)=s,  \phi_3(s,t)=s^2
$
and
$
\psi_1(s,t)=\psi_2(s,t)=1,  \psi_3(s,t)=t^2s
$, respectively.
We have that $1\not \in \{0, \infty\}=\Lambda(A(s))\cup\Lambda(B(s))$ and
 $$
\phi_{i-1}(s,t)\mid \psi_{i}(s,t)\mid \phi_{i+1}(s,t), \quad i=1, 2, 3.
$$
Therefore, there exists a matrix pencil $P(s)\in \FF[s]^{3\times 3}$  such that
$\rank P(s)\leq 1$ and  $A(s)+P(s)\se B(s)$.
In fact,
$$
A(s)+\begin{bmatrix}0&0&0\\0&0&0\\0&1&-s\end{bmatrix}=
\begin{bmatrix}s&0&0\\0&s&1\\0&1&0\end{bmatrix}\se B(s).
$$

\end{example}

\begin{theorem}\label{teoeqmult}
 Let $A(s), B(s)\in \FF[s]^{n \times n}$ be regular matrix pencils.
Let $\phi_1(s,t)\mid \dots\mid \phi_n(s, t)$ and $\psi_1(s,t)\mid \dots\mid \psi_n(s, t)$ be the homogeneous invariant factors of 
$A(s)$ and $B(s)$, respectively, and assume that  
for some $\lambda_0\in \FF\cup \{\infty\}$,
$$
m_i(\lambda_0, A(s))=m_i(\lambda_0, B(s)), \quad 1\leq i \leq n.
$$
Let $r$ be a nonnegative integer. 
There exists a matrix pencil  $P(s)\in \FF[s]^{n \times n}$ such that $\rank (P(s))\leq r$ and $A(s)+P(s)\se B(s)$ if and only if
(\ref{interlacinghomogr1}) holds.
\end{theorem}

{\bf Proof.}
If
\begin{equation} \label{emescero}
m_i(\lambda_0, A(s))=m_i(\lambda_0, B(s))=0, \quad 1\leq i \leq n,
\end{equation}
then $\lambda_0\not \in \Lambda(A(s))\cup \Lambda(B(s))$ and we can apply Theorem \ref{mainth}. In fact, condition (\ref{emescero}) implies  the condition $ \FF\cup \{\infty\} \not \subseteq \Lambda(A(s))\cup \Lambda(B(s))$ of Theorem \ref{mainth}.

In any other case,
let $\mu_a(\lambda_0, A(s))=\mu_a(\lambda_0, B(s)):=n_1$. Then
$$
A(s)\se A'(s)=\begin{bmatrix}A_{11}(s)&0\\0&A_{22}(s)\end{bmatrix}, \quad
B(s)\se B'(s)=\begin{bmatrix}A_{11}(s)&0\\0&B_{22}(s)\end{bmatrix},
$$
where $A_{11}(s)\in  \FF[s]^{n_1 \times n_1}$, $A_{22}(s), B_{22}(s)\in  \FF[s]^{(n-n_1) \times (n-n_1)}$ are regular pencils  such that $\Lambda(A_{11}(s)=\{\lambda_0\}$,
$m_i(\lambda_0, A_{11}(s))=m_i(\lambda_0, A(s))=m_i(\lambda_0, B(s))$ for $ 1\leq i \leq n_1$ and the 
homogeneous invariant factors of  $A_{22}(s)$ and $B_{22}(s)$ are
$\phi'_i(s,t)=\phi_{n_1+i}(s,t)$ and $\psi'_i(s,t)=\psi_{n_1+i}(s,t)$, $1\leq i \leq n-n_1$,  respectively.
From (\ref{interlacinghomogr1}),
$$
\phi'_{i-r}(s,t)\mid \psi'_{i}(s,t)\mid \phi'_{i+r}(s,t), \quad 1\leq i \leq n-n_1.
$$
By Theorem \ref{mainth},
there exists a matrix pencil $P_{22}(s)\in \FF[s]^{(n-n_1)\times(n-n_1) }$  such that
$\rank P_{22}(s)\leq r$ and 
$A_{22}(s)+P_{22}(s)\se B_{22}(s)$. Taking 
$P(s)= \begin{bmatrix}0&0\\0&P_{22}(s)\end{bmatrix}\in \FF[s]^{n\times n}$, we have  that
$A'(s)+P(s)\se B(s)$ and $\rank P(s)\leq r$. 

\hfill $\Box$

\begin{example}\label{secondexample}
Assume that $\FF=\ZZ_2$, $n=5$, and 
the homogeneous invariant factors of $A(s)$ and $B(s)$ are
$
\phi_1(s,t)=\phi_2(s,t)=\phi_3(s,t)=1,  \phi_4(s,t)=s(s-t),  \phi_5(s,t)=s^2(s-t)
$
and
$
\psi_1(s,t)=\psi_2(s,t)=\psi_3(s,t)=1, \psi_4(s,t)=(s-t),  \psi_5(s,t)=t^2s(s-t)
$, respectively. Then
$$
\Lambda(A(s))=\{0, 1\}, \quad \Lambda(B(s))=\{0, 1, \infty\},
$$
and  $ \FF\cup \{\infty\}= \Lambda(A(s))\cup \Lambda(B(s))=\{0,1, \infty\}$.
But
$$
(m_1(1, A(s)), \dots, m_5(1, A(s)))=
(m_1(1, B(s)), \dots, m_5(1, B(s)))=(0,0,0,1,1).
$$
Then,
$$
A(s)\se A'(s)=\begin{bmatrix}A_{11}(s)&0\\0&A_{22}(s)\end{bmatrix}, \quad
B(s)\se B'(s)=\begin{bmatrix}A_{11}(s)&0\\0&B_{22}(s)\end{bmatrix},
$$
where $A_{11}(s)=\begin{bmatrix}s-1&0\\0&s-1\end{bmatrix}$, $A_{22}(s), B_{22}(s) \in \FF[s]^{3\times 3}$ and the homogeneous invariant factors of 
$A_{22}(s)$ and $B_{22}(s)$ are
$
\phi'_1(s,t)=1,  \phi'_2(s,t)=s,  \phi'_3(s,t)=s^2
$
and
$
\psi'_1(s,t)=\psi'_2(s,t)=1,  \psi'_3(s,t)=t^2s
$, respectively.
Hence (see Example \ref{firstexample}) there exists a matrix pencil $P_{22}(s)\in \FF[s]^{3\times 3}$  such that
$\rank P_{22}(s)\leq 1$ and 
$A_{22}(s)+P_{22}(s)\se B_{22}(s)$. Taking 
$P(s)= \begin{bmatrix}0&0\\0&P_{22}(s)\end{bmatrix}\in \FF[s]^{(2+3)\times (2+3)}$, we have  that
$A'(s)+P(s)\se B(s)$ and $\rank P(s)\leq 1$. 

\end{example}

\section{Eigenvalue placement for regular matrix pencils under low rank perturbations}
\label{secplacement}

In this section we solve  Problem \ref{prdet}.

For $\lambda \in \overline{\FF}\cup \{\infty\}$, 
we will denote by $M_r(\lambda, A(s)):=\sum_{i=n-r+1}^nm_i(\lambda, A(s))$
and $M_r(A(s)):=\sum_{\lambda\in \Lambda(A(s))}M_r(\lambda, A(s))$.
\begin{proposition}\label{propnecdet}
Let $A(s), P(s)\in \FF[s]^{n \times n}$ be matrix pencils with  $ \rank (P(s))= r$. Assume that $A(s)$ and $A(s)+P(s)$ are regular. 
Let $\alpha_1(s)\mid \dots\mid \alpha_n(s)$ be the   invariant factors of $A(s)$ and  $p(s)=\det(A(s)+P(s))$. Then
\begin{equation}\label{eqdetpol}
\begin{array}{l}
\alpha_1(s)\dots\alpha_{n-r}(s)\mid p(s),\\
\mu_a(\infty, A(s))-M_r(\infty, A(s))\leq n-\deg(p(s)).
\end{array}
\end{equation}
\end{proposition}

{\bf Proof.}
Let $B(s)= A(s)+P(s)$ and let $\beta_1(s)\mid \dots\mid \beta_n(s)$ be its   invariant factors.
By
Proposition \ref{propnecesr} and Lemma \ref{lemmadiv}, condition (\ref{eqinalphar}) and
\begin{equation}
\label{eqinmr}
m_{i-r}(\infty, A(s))\leq m_{i}(\infty, B(s))\leq m_{i+r}(\infty, A(s)), \quad 1 \leq i\leq n
\end{equation}
hold.
As $p(s)=k\beta_1(s)\dots  \beta_n(s)$, with $0\neq k \in \FF$, we have that
$$
\begin{array}{l}
\alpha_1(s)\dots\alpha_{n-r}(s)\mid p(s),\\
\sum_{i=1}^{n}m_{i-r}(\infty, A(s))
\leq \sum_{i=1}^{n}m_i(\infty, B(s)).
\end{array}
$$
But
$
\sum_{i=1}^{n}m_{i-r}(\infty, A(s))=
\mu_a(\infty, A(s))-M_r(\infty, A(s))$ and $
\sum_{i=1}^{n}m_i(\infty, B(s))=
\mu_a(\infty, B(s))=
n-\deg(p(s)).
$

\hfill $\Box$

From Proposition \ref{propnecdet} we derive in Corollary \ref{cordetnecm} a lower bound for the algebraic multiplicity of the eigenvalues of the perturbed pencil.

\begin{corollary}\label{cordetnecm}
  Let
 $A(s), P(s)\in \FF[s]^{n \times n}$ be matrix pencils with $ \rank (P(s))=r$. Assume that $A(s)$ and $A(s)+P(s)$ are regular. 
 Then
 \begin{equation}\label{eq1718Getrsemi}
\mu_a(\lambda, A(s))-M_r(\lambda, A(s))\leq \mu_a(\lambda, A(s)+P(s)),
\quad \lambda \in \overline{\FF}\cup \{\infty\}.
\end{equation}
  \end{corollary}

\begin{remark}
Under the conditions of Corollary \ref{cordetnecm}, if for $\lambda \in \overline{\FF} \cup\{\infty\}$ we take
$
  d(\lambda)=\mu_a(\lambda, A(s)+P(s))-\mu_a(\lambda, A(s))+M_r(\lambda, A(s))
$ 
  and denote $\mathcal{F}=\Lambda(A(s))\cup \Lambda(A(s)+P(s))\cup \{\infty\},
$
  then
$$
\sum_{\lambda\in\mathcal{F}}d(\lambda)=
n-n+\sum_{\lambda\in\Lambda(A(s))} M_r(\lambda, A(s))=
M_r( A(s)).
$$
Condition  (\ref{eq1718Getrsemi}) implies
$d(\lambda)\geq 0$ for all $\lambda \in \overline{\FF}\cup \{\infty\}$.
Hence,
$$d(\lambda)\leq M_r(A(s)), \quad \lambda \in \overline{\FF}\cup \{\infty\}.$$
Therefore,  (\ref{eq1718Getrsemi}) implies
\begin{equation}\label{eq1718Getrsemi2}
  \mu_a(\lambda, A(s)+P(s))\leq \mu_a(\lambda, A(s))-M_r(\lambda, A(s))+M_r(A(s)),
  \quad \lambda \in \overline{\FF}\cup \{\infty\}.
\end{equation}
We notice that conditions (\ref{eq1718Getrsemi}) and (\ref{eq1718Getrsemi2}) for $r=1$ are those  of \cite[Theorem 4.4]{GeTr17}.
  \end{remark}

The next theorem gives a solution to Problem \ref{prdet} under a restriction on the field $\FF$.
We will use the following notation for a given $p(s)\in \FF[s]$ with $\deg(p(s))\leq n$,
$$\Lambda(p(s)):=\{\lambda \in \overline{\FF}\; : \; p(\lambda)=0\}\mbox{  if } \deg(p(s))=n,$$
$$\Lambda(p(s)):=\{\lambda \in \overline{\FF}\; : \; p(\lambda)=0\}\cup \{\infty\}\mbox{  if } \deg(p(s))<n.$$
\begin{theorem}\label{mainthdet}
 Let $A(s)\in \FF[s]^{n \times n}$ be a regular matrix pencil and $\alpha_1(s)\mid \dots\mid \alpha_n(s)$ be its   invariant factors. Let $p(s)\in \FF[s]$ be a nonzero monic polynomial with $\deg (p(s))\leq n$.
Assume that  $ \FF\cup \{\infty\}\not \subseteq \Lambda(A(s))\cup \Lambda(p(s))  $.
Let $r$ be  a nonnegative integer.
There exists a matrix pencil  $P(s)\in \FF[s]^{n \times n}$ such that $\rank (P(s))\leq r$ and $\det(A(s)+P(s))=kp(s)$ with $0\neq k \in \FF$ if and only if
(\ref{eqdetpol}) holds.
\end{theorem}
    {\bf Proof.}
    The necessity follows from Proposition \ref{propnecdet}.
    
    Assume now that (\ref{eqdetpol}) holds.
    Then there exists $\gamma(s)\in \FF[s]$ such that $$p(s)=\alpha_1(s)\dots\alpha_{n-r}(s)\gamma(s),$$
    and
    $$d= n-\deg(p(s))-\mu_a(\infty, A(s))+M_r(\infty, A(s))\geq 0.$$
    Let us define
    $$
\beta_i(s):=\alpha_{i-r}(s), \; 1\leq i \leq n-1, \quad \beta_n(s):=\alpha_{n-r}(s)\gamma(s),
$$
$$
m'_{i}:=m_{i-r}(\infty, A(s)), \; 1\leq i \leq n-1, \quad m'_n:=m_{n-r}(\infty, A(s))+d.
$$
Then
$$\beta_1(s) \mid \dots \mid \beta_n(s), \quad m'_1\leq \dots \leq m'_{n},$$
$$
\sum_{i=1}^n\deg(\beta_i(s))+\sum_{i=1}^n m'_i=\deg(p(s))+\mu_a(\infty, A(s))-M_r(\infty, A(s))+d=n.
$$

Let  $B(s)$ be a pencil with invariant factors $\beta_1(s) \mid \dots \mid \beta_n(s)$ and $m_i(\infty, B(s))=m'_i$, $1\leq i \leq n$.
Then,  $B(s)$ is regular and satisfies (\ref{eqinalphar}) and (\ref{eqinmr}).

By Theorem \ref{mainth}, there exists a pencil  $P(s)\in \FF[s]^{n \times n}$ such that  $\rank (P(s))\leq r$ and $A(s)+P(s)\se B(s)$.
Then there exist $0\neq k_1, k_2 \in \FF$ such that
$$
\det(A(s)+P(s))=k_1\det(B(s))=k_1k_2\beta_1(s)\dots \beta_n(s)=k_1k_2\alpha_1(s)\dots\alpha_{n-r}(s)\gamma(s)$$$$=k_1k_2p(s),\quad 0\neq k_1k_2\in \FF.
$$

\hfill $\Box$

\section{Conclusions}
\label{secconclusions}

Given a regular matrix pencil, we have completely characterized the 
Weiestrass structure of a regular pencil obtained by a low rank perturbation of it.
The conditions obtained are expressed both in terms of the homogeneous invariant factors and in terms of the partial multiplicities of the  eigenvalues.
The necessity of the conditions holds over  arbitrary fields and the sufficiency over fields with sufficient number of elements.

We also  determine the possible eigenvalues and their algebraic multiplicities of a regular matrix pencil obtained by a low rank perturbation of another regular one.
An open  problem related to this study is  to develop efficient methods for constructing the pencils $P(s)$ whose existence is stablished in Theorem \ref{mainthdet}.

For singular matrix pencils, the general case remains open. As mentioned in the Introduction section, for a particular kind  of pencils the problem is solved in \cite{DoSt14} and the generic behavior has been described in \cite{TeDo07}.

  \bibliographystyle{acm} 
\bibliography{referencesreg}

\end{document}